\documentclass[draft,12pt]{amsart}
\usepackage[all]{xy}
\usepackage{a4wide}
\usepackage{mathtools}
\usepackage{amssymb}
\usepackage{amsthm}
\usepackage{amsmath}
\usepackage{amscd}
\usepackage{verbatim}
\usepackage[all]{xy}

\usepackage{comment}

\usepackage{xfrac}
\usepackage{enumitem}

\numberwithin{equation}{section}

\usepackage{xcolor}

\newcommand{\SL}{\operatorname{SL}}
\newcommand{\Mp}{\operatorname{Mp}}
\newcommand{\rot}{\operatorname{rot}}
\newcommand{\Meas}{\operatorname{Meas}}
\newcommand{\FS}{\operatorname{FS}}

\newcommand{\tr}{\operatorname{tr}}
\newcommand{\RR}{\mathbb{R}}

\newcommand{\ZZ}{\mathbb{Z}}
\newcommand{\NN}{\mathbb{N}}
\newcommand{\HH}{\mathbb{H}}
\newcommand{\CC}{\mathbb{C}}
\newcommand{\G}{\Gamma}

\newcommand{\vol}{\operatorname{vol}}

\newcommand{\id}{\operatorname{id}}
\newcommand{\diag}{\operatorname{diag}}
\newcommand{\dd}{\ts\mathrm{d}}
\newcommand{\ts}{\hspace{0.5pt}}

\renewcommand{\hat}{\widehat}

\makeatletter
\newcommand*\defbb[1]{
	\expandafter\newcommand\csname I#1\endcsname{\mathbb{#1}}}
\newcommand*\defbbs[1]{
	\@for\@i:=#1\do{\expandafter\defbb\expandafter{\@i}}}
\makeatother
\defbbs{A,B,C,D,E,F,G,H,I,K,L,M,N,O,P,Q,R,S,T,U,V,W,X,Y,Z} 

\makeatletter
\newcommand*\deffrak[1]{
	\expandafter\newcommand\csname frak#1\endcsname{\mathfrak{#1}}}
\newcommand*\deffraks[1]{
	\@for\@i:=#1\do{\expandafter\deffrak\expandafter{\@i}}}
\makeatother
\deffraks{a,b,c,d,e,f,g,h,i,j,k,l,m,n,o,p,q,r,s,t,u,v,w,x,y,z}

\makeatletter
\newcommand*\defcal[1]{
	\expandafter\newcommand\csname cal#1\endcsname{\mathcal{#1}}}
\newcommand*\defcals[1]{
	\@for\@i:=#1\do{\expandafter\defcal\expandafter{\@i}}}
\makeatother
\defcals{A,B,C,D,E,F,G,H,I,J,K,L,M,N,O,P,Q,R,S,T,U,V,W,X,Y,Z}

\newtheorem{theorem}{Theorem}[section]
\newtheorem{lemma}[theorem]{Lemma}
\newtheorem{proposition}[theorem]{Proposition}
\newtheorem{corollary}[theorem]{Corollary}

\theoremstyle{definition}
\newtheorem{definition}[theorem]{Definition}
\newtheorem{example}[theorem]{Example}
\theoremstyle{remark}
\newtheorem{remark}[theorem]{Remark}

\usepackage[draft=false]{hyperref}

\title[Measures, modular forms, and summation formulas of Poisson type]{Measures, modular forms, \\[2mm] and summation formulas of Poisson type}

\author{Claudia Alfes}
\address{Universität Bielefeld, Fakultät für Mathematik, Postfach 100 131, 33501 Bielefeld, Germany}
\email{alfes@math.uni-bielefeld.de}

\author{Paul Kiefer}

\email{pkiefer@math.uni-bielefeld.de}

\author{Jan Maz\'{a}\v{c}}

\email{jmazac@math.uni-bielefeld.de}

\thanks{The authors were funded by the Deutsche Forschungsgemeinschaft (DFG, German Research Foundation) – Project-ID 491392403 – TRR 358.}

\date{\today}

\begin{document}

\begin{abstract}
    In this article, we show that Fourier eigenmeasures supported on spheres with radii given by a locally finite sequence, which we call $k$-spherical measures, correspond to Fourier series exhibiting a modular-type transformation behaviour with respect to the metaplectic group. A familiar subset of such Fourier series comprises holomorphic modular forms. This allows us to construct $k$-spherical eigenmeasures and derive Poisson-type summation formulas, thereby recovering formulas of a similar nature established by Cohn--Gonçalves, Lev--Reti, and Meyer, among others. Additionally, we extend our results to higher dimensions, where Hilbert modular forms yield higher-dimensional $k$-spherical measures.
\end{abstract}

\maketitle

\section{Introduction}

Our work connects the theory of modular forms to mathematical physics, namely aperiodic order and crystallography, and to harmonic analysis. 
We focus on two aspects: (1) The characterization of certain eigenmeasures of the Fourier transform in terms of modular forms, and (2) Poisson-type summation formulas, i.e., Poisson summation formulas with weights involving the coefficients of modular forms.
 
We start by describing the historical background, and the relation between the classical \emph{Poisson summation formula}, the \emph{lattice Dirac comb} measure, and the \emph{Jacobi theta function}. We illustrate our results by examining a special tempered distribution initially studied by Guinand \cite{guinand}. We close the introduction by stating our main results.

\subsection{Background: Eigenmeasures of the Fourier transform}
The classical Fourier transform is a unitary operator of order $4$ acting on $L^{2}(\RR)$ with eigenvalues $\{ \pm 1, \, \pm i\}$. Its eigenfunctions in 1D, the Hermite functions, are well-studied.

The discovery of quasicrystalline materials in the 80s, namely non-periodic atomic structures with discrete diffraction patterns, led to a variety of problems asking for generalizations of these results beyond the setting of finite, absolutely continuous measures, see \cite{meyer1995,lagarias}, or the monographs \cite{TAO1,TAO2} for general background.

A famous extension to the realm of unbounded, translation-bounded measures is given by the \emph{lattice Dirac comb} 
\[
\delta^{}_{\ZZ} \, = \, \sum^{}_{m\in \ZZ} \delta^{}_{m},
\]
where $\delta^{}_{x}$ is the Dirac point measure centred at $x$.
This measure satisfies $\widehat{\delta}_{\ZZ} \, = \, \delta^{}_{\ZZ}$ in the sense of tempered distributions.
This follows immediately from the Poisson summation formula 
\[
 \sum_{m\in\ZZ}\, \varphi(m) \ = \ \sum_{m\in\ZZ}\, \hat{\varphi}(m)\, ,
\]
which holds for all Schwartz functions $\varphi \in\calS(\RR)$. The Poisson summation formula $\widehat{\delta}_{\ZZ} \, = \, \delta^{}_{\ZZ}$ also holds in the stronger sense of measures \cite{FAunbouned}.

This measure belongs to the class of \emph{crystalline measures} \cite{meyerpnas}, which are one-dimensional, tempered measures such that their support and the support of their Fourier transform are both sparse. Such measures are also called \emph{doubly-sparse measures} \cite{baakeeigenmeasure} and include the so-called \emph{Fourier quasicrystals} in the sense of \cite{kurasovsarnak}. 

Lev and Olevskii \cite{levolevskii} characterised crystalline measures of $\RR$ and positive measures on $\RR^n$ whose support and spectrum are uniformly discrete. Recall that a set is called uniformly discrete if the distance between two points is bounded from below by a positive constant. They show that a crystalline measure supported on a uniformly discrete set can be obtained from the Poisson summation formula using dilation and a finite number of shifts, modulations, and by taking suitable linear combinations.

Recently, Baake, Spindeler and Strungaru classified all one-dimensional periodic Fourier eigenmeasures and those with uniformly discrete support \cite{baakeeigenmeasure}. A further interesting example of an unbounded measure that is a Fourier eigenmeasure, in fact, an eigendistribution, goes back to Guinand \cite{guinand}. Meyer  \cite[Theorem~4.1]{meyerpnas} used Guinand's work to construct a Fourier eigenmeasure with eigenvalue $-i$.

On the other hand, Lev and Olevskii \cite{LO_fourier_17} and Strungaru \cite{Strungaru_weighted_DC} showed that the support of a crystalline measure cannot be a~model set \cite[Def.~7.2]{TAO1} if one goes beyond the generalized Dirac combs (in the sense of \cite[Def.~1.2]{Meyer_almostperiodic}). Further examples of such measures can be found in  \cite{LO_crystalline_13, LO_qc_discrete16,Kolountzakis,Favorov}, or \cite{meyerpm} for a new two-dimensional extension. A full classification of these measures is still lacking, but the first steps have already been taken in \cite{baakepurepoint}.

In this paper, we give a novel perspective on the above story. We show that modular forms, well-known objects from number theory, give rise to eigenmeasures under the Fourier transform.

\subsection{Background: Jacobi's theta function and Poisson's summation formula}
The Poisson summation formula implies the modular transformation behaviour of the \emph{Jacobi theta function}
\begin{equation}\label{eq:theta}
 \vartheta(\tau)\ = \ \sum_{m\in\ZZ}\ts e( m^2 \tau)\, ,\quad \mbox{with } \ \tau=u+i\ts v \in \HH.
\end{equation}
Here, $\HH$ denotes the upper half of the complex plane, and from now on, we abbreviate $e\left(y\right)=e^{2\pi i y}$. The Jacobi theta function is an example from a broad class of functions called \emph{modular forms}. 
The latter are holomorphic functions on $\HH$ that transform nicely with respect to the fractional linear transformations induced by the action of discrete subgroups of $\SL_2(\RR)$ on $\HH$ (compare Definition \ref{def:modformthetamult}). 

Applying the Poisson summation formula to the Schwartz function $f_a(x)=e^{-\pi x^2/a}$ with $a>0$ gives
\[
\sum_{m\in\ZZ}\ts e^{-\frac{\pi m^2}{a}} \, = \, \sqrt{a} \sum_{m\in\ZZ}\ts  e^{-a\pi m^2}\, ,
\]
since $\hat{f}_a(y)=\sqrt{a}\ts f_{1/a}(y)$. This implies
\begin{equation}
    \label{eq:Jacobi_transformation}
 \vartheta\!\left({-}\frac{1}{4\tau}\right) \, = \, \sqrt{-2i \tau} \, \vartheta(\tau) \ts
\end{equation}
for $\tau = ia/2$ with $a>0$. The general case follows from the identity theorem.

This functional equation gives the modular transformation behaviour of $\vartheta(\tau)$ for the matrix $\left(\begin{smallmatrix}
    0& -1/2\\2&0
\end{smallmatrix}\right)$. Moreover, it is $1$-periodic, holomorphic on $\HH$ and in the cusps $\infty,\,0,\,1/2$. In other words, $\vartheta$ is a holomorphic elliptic modular form of weight $1/2$ for a certain discrete subgroup of $\SL_2(\RR)$.

\subsection{Background: Poisson-type summation formulas and crystallography}
In \cite{meyerpnas}, Meyer showed that there are variations of the Poisson summation formula. In particular, he studied the summation formula arising from Guinand's distribution. It involves as coefficients the representation number of $3$ squares (compare \eqref{eq:psfguinand}).

Further Poisson-type summation formulas appeared in the context of harmonic analysis, for instance, in the work of  Cohn--Gon\c{c}alves \cite{cohngonzalves}, Gon\c{c}alves \cite{goncalves}, or Lev--Reti \cite{levreti}. A joint feature of these formulas is that they involve coefficients of modular forms.
Our results give a common framework for such Poisson-type summation formulas.

From this point of view, our work deepens the connection between measures of crystalline type and modular forms. This connection was initiated by the seminal work of Viazovska \cite{maryna1,maryna2} and has led to intriguing results touching upon the connections between measures, Poisson summation formulas, Fourier interpolation formulas, Fourier uniqueness, and modular forms \cite{radvia,cohngonzalves,stoller,RadchenkoStoller,goncalves,gerbellivenkatesh}. 

We mention that the Poisson summation formula also plays a crucial role in the theory of lattice combs and crystallography, as it brings the key to the pure point diffraction measure of perfect crystals, whose atomic configurations are described via lattices. The theory of aperiodic order uses a projection of higher-dimensional lattices to describe mathematical quasicrystals \cite[Sec.~7]{TAO1}. In \cite{PSF_quasicrystals_Nicu_Christoph} and \cite{modulated_crystals}, it was shown that Poisson summation formulas of the underlying higher-dimensional lattices and the pure-point diffraction formula for the quasicrystals (regular model sets) can be derived from one another.

\subsection{Guinand's tempered distribution and the third power of Jacobi's theta function}\label{subsec:guinand}
We illustrate the novel viewpoint of our work for a tempered distribution first considered by Guinand \cite{guinand}, which has eigenvalue $-i$ under the Fourier transform. In Theorem \ref{lem:IsoBetweenMeasuresAndFunctions}, we give an isomorphism between a certain space of measures and Fourier series satisfying a modular-type functional equation. In particular, this gives a way to construct Guinand's tempered distribution \cite{guinand} via the third power of the Jacobi theta function
\[
\label{eq:Jacobi_cube}
\vartheta^3(\tau)\,= \, 1+\sum_{m=1}^\infty r_3(m) \ts e(m\tau),
\]
where $r_3(m)=\left|\left\{(x_1,\, x_2,\, x_3)\in\ZZ^3\,|\, x_1^2+x_2^2+x_3^2=m\right\}\right|$ denotes the representation number of $m$ as a sum of $3$ squares. 
Consider $\vartheta^3(\tau / 2)$ and define the measure $\widetilde{\sigma}_{3}$ in $\RR^3$ with coefficients $r_3(m)$ as 
\[\widetilde{\sigma}_{3} \, = \, \delta^{}_{0} + \sum_{m=1}^{\infty} \, \frac{r^{}_{3}(m)}{\vol(S^{3}_m)} \, \delta^{}_{S^{3}_m},\]
where $\delta^{}_{S^{3}_m}$ is the surface measure that integrates over the sphere $S^{3}_m$ of radius $\sqrt{m}$. The modular transformation behaviour of $\vartheta(\tau)$ implies 
\[ \vartheta^3\!\left( {-}\frac{1}{4\tau} \right) \ = \ \sqrt{-2i\tau\ts}^{\ts 3} \vartheta^3 \! \left(\tau \right), \]
which encodes, again by Theorem \ref{lem:IsoBetweenMeasuresAndFunctions}, that $\widetilde{\sigma}_{3}$ is an eigenmeasure under the Fourier transform, still with eigenvalue 1.

Moreover, by considering the radial test functions and the descent methods for them (see Section \ref{sec:MeasuresFromFourierSeries}), we obtain the radial part 
\[
 \sigma_3 \, = \, -\delta_0'+\sum_{m=1}^\infty \, \frac{r_3(m)}{\sqrt{m}} \,\delta^{}_{\sqrt{m}},
\]
which is Guinand's distribution \cite{guinand}. We recover the following summation formula of Poisson type (see \ref{lem:ModularFormsAnd1DimOddMeasures} and \ref{rem:thetak}, also compare the work of Lev and Reti \cite{levreti})
\begin{align}\label{eq:psfguinand}
 \varphi'(0)+\sum_{m=1}^\infty \frac{r_3(m)}{\sqrt{m}}\,\varphi(\sqrt{m}\,) \, = \,  i\left(\hat\varphi'(0) +\sum_{m=1}^\infty \frac{r_3(m)}{\sqrt{m}}\,\hat\varphi(\sqrt{m}\,)\right).
\end{align}

\subsection{Main results}
We now turn to the statement of our main results.
For $\lambda \in \IR_{\geq 0}$, we let $S_{\lambda}^k \subseteq \IR^k$ be the sphere with radius $\sqrt{\lambda}$. For $\lambda = 0$, we let $S_{\lambda}^k = \{0\}$. We write $\delta_{S_\lambda^k}$ for the surface measure that integrates over the sphere $S_\lambda^k$. We recall that a sequence is locally finite if it has no finite accumulation points.

\begin{definition}
Let $\Meas_k$ be the space of measures $\mu$ on $\IR^k$ such that there exist locally finite sequences $\lambda_m,\, \hat{\lambda}_m \in \IR_{\geq 0}$ and at most polynomially growing sequences $c_m, \hat{c}_m \in \IC$ with
\[\mu = \sum_{m = 0}^\infty c_m \delta_{S_{\lambda_m}^k} \qquad \text{and} \qquad \hat{\mu} = \sum_{m = 0}^\infty \hat{c}_{m} \delta_{S_{\hat{\lambda}_m}^k}.\]
We call these measures \emph{$k$-spherical measures}.
\end{definition}
In particular, $\Meas_1$ is the space of all crystalline measures on $\RR$. 

\begin{definition}
Let $\FS_k$ be the set of holomorphic functions $f : \IH \to \IC$ such that there exist locally finite sequences $\lambda_m, \hat{\lambda}_m \in \IR_{\geq 0}$ and at most polynomially growing sequences $c_m, \hat{c}_m \in \IC$ with
\[f(\tau) = \sum_{m = 0}^\infty c_m e\left(\frac{\lambda_m \tau}{2} \right)\]
such that $\tilde{f}(\tau) = \sqrt{-i \tau}^{\, -k} f(-1/\tau)$ has a (convergent) Fourier expansion
\[\tilde{f}(\tau) = \sum_{m = 0}^\infty \hat{c}_m e\bigg(\ts \frac{\hat{\lambda}_m \tau }{2}\bigg).\]
\end{definition}
\begin{remark}
If $f=\tilde{f}$ in the above definition, $f$ is a holomorphic function on the upper half plane $\HH$ and in the cusp $\infty$, and satisfies a modular transformation behaviour for the matrix $\left(\begin{smallmatrix} 0&-1\\1&0\end{smallmatrix}\right)\in \SL_2(\ZZ)$ (compare \eqref{eq:slashoperator} and Definition \ref{def:modformthetamult}). 
\end{remark}

Let $\Mp_{2}(\IR)$ be the metaplectic cover of $\SL_{2}(\IR)$, i.e.,\@ elements in $\Mp_2(\IR)$ are pairs $(M, \phi)$ with $M = \left(\begin{smallmatrix} a & b \\ c & d\end{smallmatrix}\right) \in \SL_2(\IR)$ and $\phi : \IH \to \IC$ is a holomorphic function satisfying $\phi^2(\tau) = c \tau + d$. The multiplication is given by
\[(M_1, \phi_1) (M_2, \phi_2) \, = \, \bigl(M_1 M_2,\, \phi_1(M_2(\cdot)) \phi_2 \bigr).\]

\noindent For a tempered distribution $\mu \in \calS'(\IR^k)$ and $M \in \Mp_2(\IR)$, we define $\mu \vert M$ by
\begin{equation}\label{eq:muM}
    \langle \mu \vert M, \, \varphi \rangle \, = \,  \langle \mu, \, \omega(M^{-1}) \varphi \rangle,
    \end{equation}
where $\omega$ is the Schr\"odinger model of the Weil representation which gives the action of $\Mp_2(\IR)$ on Schwartz functions (compare \eqref{eq:WeilRep}). Moreover, let its theta function be defined by
\begin{equation}\label{eq:thetamu}
    \theta_\mu(\tau) \, = \, \langle \mu, \, g(\cdot, \tau) \rangle,
    \end{equation}
where $g(x, \tau) = e\bigl(x^{\top}\!x \ts \tau / 2\bigr)$ is the Gaussian. The group $\Mp_2(\RR)$ acts on $\theta_\mu(\tau)$ by the usual slash operator $\vert_{k/2}$ of weight $k / 2$ (compare \eqref{eq:slashoperator}).

\begin{theorem}\label{lem:IsoBetweenMeasuresAndFunctions}
The map
\[\theta : \Meas_k \, \to \, \FS_k, \quad \mu \, \mapsto \, \theta_\mu\]
is an $\Mp_2(\IR)$-equivariant isomorphism, i.e., $\theta_{\mu \vert M} = \theta_\mu \vert_{k / 2} M$ for all $M \in \Mp_2(\IR)$. The inverse is given by
\[f(\tau) = \sum_{m = 0}^\infty c_m e\left(\frac{\lambda_m \tau}{2}\right) \ \longmapsto \ \mu_f = \sum_{m = 0}^\infty \frac{c_m}{\vol(S_{\lambda_m}^k)} \delta_{S_{\lambda_m}^k}.\]
\end{theorem}

\begin{remark}
     Under the above isomorphism, measures $\mu \in \Meas_k$ with $\hat{\mu} = \mu$ correspond to functions $f \in \FS_k$ with $f(-1/\tau) = \sqrt{-i\tau}^{\ts k} f(\tau)$. As there are infinitely many modular forms satisfying this property, we can construct infinitely many measures of this type. 
  
  More generally, measures $\mu \in \Meas_k$ with $\hat{\mu} = - \mu$ or $\hat{\mu} = \pm i \mu$ can be constructed using modular forms transforming with respect to the Fricke group with a more general multiplier system (compare Section \ref{sec:poissontype}).
\end{remark}

\begin{remark}
    In Section \ref{sec:poissontype} we combine Theorem \ref{lem:IsoBetweenMeasuresAndFunctions} with results of Lev and Reti \cite{levreti} to obtain formulas of Poisson type that involve Fourier coefficients of modular forms as weights.
\end{remark}

\begin{remark}
    In Section \ref{sec:measuresfourierhigher}, we generalize Theorem \ref{lem:IsoBetweenMeasuresAndFunctions} to higher dimensions, i.e.\@, to measures in the space $\Meas_{\underline{k}}$, where $\underline{k}=(k_1,\ldots,k_n)$. We also give a relation to Hilbert modular forms.
\end{remark}
We illustrate Theorem \ref{lem:IsoBetweenMeasuresAndFunctions} with two examples.
\begin{example}
\begin{enumerate}
    \item Consider Dedekind's Delta function $\Delta(\tau)$ given by
    \[
    \Delta(\tau)\, = \, e(m\tau)\prod_{m=1}^\infty \bigr(1-e(m\tau)\bigl)^{24}\, =\, \sum_{m=1}^\infty \tau(m)\ts e(m\tau),
    \]
    where we denote the Fourier coefficients of $\Delta$ by $\tau(m)$ due to historical reasons. It is a cusp form of weight $12$ and satisfies $\Delta(-1/\tau)=\tau^{12}\Delta(\tau)$. By Theorem \ref{lem:IsoBetweenMeasuresAndFunctions} the measure
    \[
    \mu_\Delta \, =\, \sum_{m=0}^\infty \tau(m) \ts \delta^{}_{S^{24}_{2m}}
    \]
    satisfies $\mu_\Delta=\hat\mu_\Delta$.
    \item For the square of Jacobi's theta function
\[
    \vartheta(\tau)^2=1+\sum_{m=1}^\infty r_2(m)\ts e(m\tau),\, \text{where } r_2(m)=\lvert\{(x_1,x_2)\in\ZZ^2\,|\, x_1^2+x_2^2=m\}\rvert,
\]
we have $\vartheta^2(-1/4\tau)=\sqrt{-2i\tau}^2 \vartheta^2(\tau) $. Let $\tilde{\vartheta}(\tau) = \vartheta(\tau / 2)$. Then $\tilde{\vartheta}^2(-1/\tau) = \sqrt{-i \tau}^{\, 2} \tilde{\vartheta}^2(\tau)$.
By Theorem \ref{lem:IsoBetweenMeasuresAndFunctions}, the measure
$\mu = \sum_{m=0}^\infty r_2(m) \ts \delta^{}_{S^2_{m}}$
satisfies $\mu=\hat\mu$. This construction generalizes to powers $k\geq 2 $ (compare Remark \ref{rem:thetak}) and gives a new proof of the radial analogue of Poisson summation formula derived in \cite{baakeradial}.
\end{enumerate}
\end{example}

\begin{remark}
There is another interesting connection between crystalline measures and the theory of modular forms. Kurasov and Sarnak \cite{kurasovsarnak} construct positive crystalline measures via pairs of stable polynomials. Period polynomials of modular forms give rise to such polynomials (for the theory of period polynomials, compare \cite{cohenstroemberg}).
\end{remark}

\subsubsection{Outline}
In Section \ref{sec:MeasuresFromFourierSeries}, we recall the action of the metaplectic group $\Mp^{}_{2}(\RR)$ on Schwartz functions and how to obtain  Schwartz functions in low dimensions from those in higher, and prove the main result. In Section \ref{sec:poissontype}, we introduce modular forms and derive Poisson-type summation formulas. In particular, we recover results of this type from \cite{guinand,cohngonzalves,goncalves,levreti,meyerpnas}.  Section \ref{sec:measuresfourierhigher} is devoted to generalizing the previously derived results to higher dimensions. In the last section, we show how the usual operations on the function spaces $\FS_{\underline{k}}$ induce operations on $\Meas_{\underline{k}}$.


\section{Measures and Fourier series}\label{sec:MeasuresFromFourierSeries}

\subsection{Schwartz functions and tempered distributions}
By $\calS(\IR^k)$ we denote the space of \emph{Schwartz functions} consisting of smooth functions $\varphi$ on $\RR^k$ such that, for each $n\geq 0$ and each $m=(m_1,\ldots,m_k)$, the norm
\[
\lvert\lvert \varphi\rvert\rvert_{n,m} \, := \,  \sup_{x\in \RR^k} |x|^n\ts |\partial^m\varphi(x)|
\]
is finite. 
By $\calS'(\IR^k)$ we denote the space of \emph{tempered distributions}.

If $\mu \in \calS'(\IR^k)$ and $\varphi \in \calS(\IR^k)$, we write $\langle \mu, \varphi \rangle$ for the evaluation of $\mu$ at $\varphi$. For $\varphi \in \calS(\IR^k)$, we denote by 
\[
\hat{\varphi}(y) \, =\, \int_{\RR^k} \varphi(x) \, e(-x^{\top}\! y) \dd x
\]
its Fourier transform. For $\mu \in \calS'(\IR^k)$ we define $\hat{\mu}$ via $\langle \hat{\mu}, \varphi \rangle = \langle \mu, \hat{\varphi} \rangle$. Moreover, we define $\frac{\dd}{\dd x_i} \mu$ by $\langle \frac{\dd}{\dd x_i} \mu,  \varphi \rangle = - \langle \mu,  \frac{\dd}{\dd x_i} \varphi \rangle$.

By $O(k)$, we denote the orthogonal group of $k\times k$-matrices with real entries. We say that a Schwartz function $\varphi \in \calS(\IR^k)$ is \emph{radial} if $\varphi(g x) = \varphi(x)$ for all $x \in \IR^k,\, g \in O(k)$, and denote the space of radial Schwartz functions by $\calS(\IR^k)^{O(k)}$. Similarly, we define radial distributions and denote the space of radial distributions by $\calS'(\IR^k)^{O(k)}$.

We say that a distribution $\mu$ on $\RR$ is \emph{even} (resp.\@ \emph{odd}) if $\langle\mu,\varphi\rangle=0$ for every odd (resp.\@ even) $\varphi\in\calS(\RR)$.

\subsection{The metaplectic group and its action on Schwartz functions and Fourier series}

For $\varphi \in \calS(\IR^k)$, the Schr\"odinger model $\omega$ of the Weil representation is given by
\begin{align}
&\omega(\rot(a)) \ts \varphi(x) = \sqrt{a}^{\, -k} \varphi(x / a), \qquad \rot(a) = \left(\begin{pmatrix} a & 0 \\ 0 & a^{-1} \end{pmatrix}\!, \, a^{-1/2}\right), \nonumber\\[2pt]
&\omega(t(b)) \ts \varphi(x) = e(-(x^{\top}\! x)\ts b / 2) \ts \varphi(x), \qquad t(b) = \left(\begin{pmatrix} 1 & b \\ 0 & 1 \end{pmatrix}\!, \,  1\right), \label{eq:WeilRep}\\[2pt]
&\omega(S) \ts \varphi(x) = \hat{\varphi}(-x) = \sqrt{i}^{\, k} \int_{\IR^k} \varphi(y) \ts e(-x^{\top}\! y) \dd y, \qquad S = \left(\begin{pmatrix} 0 & -1 \\ 1 & 0 \end{pmatrix}\!, \, \sqrt{\tau}\right). \nonumber
\end{align}
Recall that the Gaussian $g(x, \tau) = e\bigl((x^{\top}\! x)\ts \tau / 2\bigr)$ is a complex-valued Schwartz function. The action of $\Mp_2(\RR)$ on the complex variable can be described via the slash operator. For $k\in \ZZ$ and $f:\HH\to\CC$ we let
\begin{equation}\label{eq:slashoperator}
 f|_{k/2} (M,\phi) = \phi(\tau)^{-k} f(M\tau),
\end{equation}
where $M\tau=\frac{a\tau+b}{c\tau+d}$ for $M=\left(\begin{smallmatrix}a&b\\c&d\end{smallmatrix}\right)$.

\begin{lemma}
For $M \in \Mp_2(\IR)$, we have $\omega(M,\phi) g(\cdot, \cdot) = g(\cdot, \cdot) \vert_{k / 2} (M,\phi)^{-1}$.
\end{lemma}

\begin{proof}
It suffices to show this for the generators. We have
\[\rot(a)^{-1} = \rot(a^{-1}), \quad t(b)^{-1} = t(-b), \quad S^{-1} = \left(\begin{pmatrix}0 & 1 \\ -1 & 0\end{pmatrix}\!, \, \sqrt{-\tau}\right)\!, \]
and thus
\begin{align*}
&\omega(\rot(a)) g(x, \tau) = \sqrt{a}^{\,-k} g(x / a, \tau) = \sqrt{a}^{\, -k} g(x, \tau / a^2) = g \vert_{k / 2} \ts \rot(a^{-1}), \\[2pt]
&\omega(t(b)) g(x, \tau) = e\bigl(-(x^{\top}|! x)\ts b / 2 \bigr)\ts g(x, \tau) = g(x, \tau - b) = g \vert_{k/2} \ts t(-b), \\[2pt]
&\omega(S) g(x, \tau) = \sqrt{i}^{\, k} \ts \hat{g}(-x, \tau) = \sqrt{i}^{\, k} \sqrt{-i \tau}^{\, -k}  g(x, -1/\tau) = g \vert_{k / 2} \ts S^{-1}.\qedhere
\end{align*}
\end{proof}

\subsection{Proof of Theorem \ref{lem:IsoBetweenMeasuresAndFunctions}}
Let $\mu \in \calS'(\IR^k)$ be a tempered distribution and consider $M \in \Mp_2(\IR)$. Recall the definition of $\mu \vert M$ and its theta function $\theta_\mu(\tau)$ from the introduction (Equations \eqref{eq:muM} and \eqref{eq:thetamu}).

\begin{lemma}\label{lem:Mp2Equivariance}
For $\mu \in \calS'(\IR^k)$ and $M \in \Mp_2(\IR)$ we have $\theta_{\mu \vert M} = \theta_\mu \vert_{k/2} M$, i.e.\@, the map $\mu \mapsto \theta_\mu$ is $\Mp_2(\IR)$-equivariant. Moreover, if $\mu \in \calS'(\IR^k)^{O(k)}$, then $\hat{\mu} = \mu$ if and only if
$$\theta_\mu(-1/\tau) \, = \, \sqrt{-i\tau}^{\, k} \theta_\mu(\tau).$$
\end{lemma}

\begin{proof}
For the first assertion, observe
\[\theta_{\mu \vert M}(\tau) \, = \, \langle \mu \vert M, \, g(\cdot, \tau) \rangle \, = \, \langle \mu, \, \omega(M^{-1}) g(\cdot, \tau) \rangle \, = \, \langle \mu, \, g \vert_{k / 2} M \rangle \, = \, \theta_\mu \vert_{k/2} M.\]
For the second claim, remark that the Gaussians $g(\cdot, \tau)$ span a dense subset of $\calS(\IR^k)^{O(k)}$ (by \cite[Lemma 2.2]{viaoptimality}, \cite[Proposition 4.1]{RadchenkoStoller}). Thus, $\hat{\mu} = \mu$ if and only if $\hat{\mu} = \mu$ on all Gaussians $g(\cdot, \tau)$. Since
\[\theta_\mu(-1/\tau) \, = \,  \sqrt{\tau}^{\,k} (\theta_\mu \vert_{k/2} \ts S)(\tau) \, = \, \sqrt{\tau}^{\, k} \theta_{\omega(S)\mu}(\tau) \, = \, \sqrt{-i\tau}^{\, k} \theta_{\ts\hat{\mu}}\ts(\tau),\]
we see that $\hat{\mu} = \mu$ if and only if $\theta_\mu(-1/\tau) = \sqrt{-i \tau}^{\, k} \theta_\mu(\tau)$.
\end{proof}

We now turn to the proof of the main result of the Introduction.

\begin{proof}[Proof of Theorem \ref{lem:IsoBetweenMeasuresAndFunctions}]
Lemma \ref{lem:Mp2Equivariance} gives the $\Mp_2(\IR)$-equivariance of the map $f \mapsto \mu_f$. Its well-definedness follows from the denseness of the span of all Gaussians. Thus, $\theta$ is an isomorphism. The correspondence between eigendistributions and functions possessing a~certain transformation formula follows again from Lemma \ref{lem:Mp2Equivariance}.
\end{proof}

\subsection{Tempered distributions in odd dimensions}\label{subsec:tempered}
By Lemma 3.2 of \cite{levreti}, the map 
\[
i_* : \calS(\IR) \to \calS(\IR^k)^{O(k)},\quad (i_* \varphi)(x) = \varphi(\lvert x \rvert)
\]
gives a continuous linear map from even functions in $\calS(\RR)$ to radial Schwartz functions on $\RR^k$. Their Fourier transform can be derived from the original even function.

\begin{lemma}[{\cite[Corollary 1.2]{GrafakosFourierRadial}\cite[Theorem 4.1 and Remark 4.3]{levreti}}]\label{lem:FourierTransformRadial}
Let $k$ be odd. Then the Fourier transform of $i_* \varphi$ is given by
\[\widehat{i_*\varphi \ts}(x) \, = \,  \sum_{j = 0}^{\frac{k - 1}{2}} \gamma_{j, k}(x) \, i_*\!\left(D_{j,k}(x) \hat{\varphi}\ts\right)(x),\]
where
\[D_{j,k}(x) = \begin{cases}\frac{\dd^j}{\dd x^j}, & \mbox{ if } x \neq 0, \\[2pt]
	\frac{\dd^{k - 1}}{\dd x^{k - 1}}, & \mbox{ if } x = 0, \, j = 0, \\[2pt]
	0, & \mbox{ otherwise},
\end{cases}\]
and
\[\gamma_{j, k}(x) = \begin{cases}\beta_{j, k} \lvert x \rvert^{j - k + 1}, & \mbox{ if } x \neq 0, \\
	\alpha_k, & \mbox{ if } x = 0,\, j = 0, \\
	0, & \mbox{ otherwise}.
\end{cases}\]
where
$$\beta_{j, k} = (-1)^{j} \frac{(k - j - 2)!}{(j - 1)! (k - 2j - 1)!!} \frac{1}{(2 \pi)^{\frac{k - 1}{2}}},\qquad n!!:= n(n-2)(n-4)\cdots,$$
and
$$\alpha_k = \frac{(-1)^{\frac{k - 1}{2}}}{(k - 2)!!} \frac{1}{(2 \pi)^{\frac{k - 1}{2}}}$$
with the convention $\beta_{0, k} = 0$ for $k \geq 3, \, \beta_{0, 1} = 1$ and $\alpha_1 = 1$.
\end{lemma}

If $\mu$ is an $O(k)$-invariant distribution on $\calS(\IR^k)$, we define an even tempered distribution $i^* \mu$ on $\calS(\IR)^{O(1)}$ via $\langle i^* \mu, \varphi \rangle = \langle \mu, i_* \varphi \rangle$. The Fourier transform of the $O(k)$-invariant test functions gives the Fourier transform of the even tempered distributions. 

\begin{lemma}\label{lem:EvenMeasureFromSphericalMeasure}
Let $k$ be odd and $\mu \in \calS'(\IR^k)^{O(k)}$. Then
$$\widehat{i^* \mu\ts } \, = \, \sum_{j = 0}^{\frac{k - 1}{2}} (-1)^j D_{j,k}(x) \ts \gamma_{j, k}(x) i^* \hat{\mu}.$$
\end{lemma}

\begin{proof}
For $\varphi \in \calS(\IR)^{O(1)}$, we have
\begin{align*}
\langle \widehat{i^* \mu}, \varphi \rangle
\, = \, \langle i^* \mu, \hat{\varphi} \,\rangle \, = \, \langle \mu, i_* \hat{\varphi} \, \rangle \, = \, \langle \ts \hat{\mu}, \widehat{i_* \hat{\varphi}\,} \,\rangle 
\, = \, \langle \ts \hat{\mu}, \sum_{j = 0}^{\frac{k - 1}{2}} \gamma_{j, k}(x) \ts i_* (D_{j,k}(x) \hat{\hat{\varphi}\,}) \rangle 
\end{align*}
by Lemma \ref{lem:FourierTransformRadial}. Therefore, we obtain
\[
 \langle \ts \sum_{j = 0}^{\frac{k - 1}{2}} (-1)^j D_{j,k}(x)\ts \gamma_{j, k}(x) \ts i^* \hat{\mu}, \, \varphi \rangle,
 \]
which shows the assertion.
\end{proof}

\begin{lemma}\label{lem:OddMeasureFromSphericalMeasure}
Let $k \geq 3$ be odd and $\mu \in \calS'(\IR^k)^{O(k)}$. Then
\[\widehat{\frac{1}{x} i^* \mu \,} \, = \, 2 \pi i\sum_{j = 0}^{\frac{k - 1}{2}} (-1)^{j+1} \tilde{D}_{j,k}(x) \ts \gamma_{j, k}(x) \ts i^* \hat{\mu},\]
where
\[\tilde{D}_{j, k}(x) = \begin{cases}
    \frac{\dd ^{j - 1}}{\dd x^{j - 1}}, & \mbox{ if } x \neq 0, \\[2pt]
	\frac{\dd^{k - 2}}{\dd x^{k - 2}}, & \mbox{ if } x = 0, j = 0, \\[2pt]
	0, & \mbox{ otherwise}.
\end{cases}\]
\end{lemma}

\begin{proof}
    We let $\varphi \in \calS(\IR)$ be odd, and obtain
    \[
    \langle \ts \widehat{\frac{1}{x} i^* \mu \,}, \varphi \rangle
        \, = \, \langle i^* \mu, \frac{\hat{\varphi}}{x} \ts \rangle 
        \, = \, \langle \ts \sum_{j = 0}^{\frac{k - 1}{2}} (-1)^j D_{j,k}(x) \ts \gamma_{j, k}(x) \ts i^* \hat{\mu},\, \widehat{\bigg(\frac{\hat{\varphi}}{x}\bigg)} \rangle.
    \]
    Since $k \geq 3$, we have $D_{j, k}(x) = \frac{\dd}{\dd x} \tilde{D}_{j, k}(x)$.  Using $\frac{\dd}{\dd x} \widehat{\left(\frac{\hat{\varphi}(x)}{x}\right)} = (-2 \pi i) \hat{\hat{\varphi}\,}$ the above expression equals
    \begin{align*}
        \langle \ts \sum_{j = 0}^{\frac{k - 1}{2}} (-1)^{j+1} \tilde{D}_{j,k}(x) \ts \gamma_{j, k}(x) \ts i^* \hat{\mu}, \, \frac{\dd}{\dd x} \widehat{\left(\frac{\hat{\varphi}}{x}\right)} \rangle \,&= \, \langle \ts \sum_{j = 0}^{\frac{k - 1}{2}} (-1)^{j+1} \tilde{D}_{j,k}(x) \ts \gamma_{j, k}(x) \ts i^* \hat{\mu}, \, (-2 \pi i)\hat{\hat{\varphi}\,} \rangle \\
        &= \langle 2 \pi i\sum_{j = 0}^{\frac{k - 1}{2}} (-1)^{j+1} \tilde{D}_{j,k}(x) \ts \gamma_{j, k}(x) \ts i^* \hat{\mu}, \, \varphi \rangle.\qedhere
    \end{align*}
\end{proof}

We restrict this result to the space $\Meas_k$ and obtain the following corollary, which allows us to construct eigenmeasures with eigenvalue $i^{\varepsilon+1}$ provided we have an eigenmeasure with eigenvalue $i^{\varepsilon}$.  

\begin{corollary}\label{cor:OddMeasureFromSphericalMeasure}
    Let $\mu\in \Meas_k$ with
    \[\mu = \sum_{m = 0}^\infty c_m \delta_{S_{\lambda_m}^k} \qquad \text{and} \qquad \hat{\mu} = \sum_{m = 0}^\infty \hat{c}_{m} \delta_{S_{\hat{\lambda}_m}^k}\]
    and assume that the radius of the sphere is $\lambda_0 = \hat{\lambda}_0 = 0$. Then the tempered distribution
    \[\nu \, = \, \frac{1}{x} i^* \mu = -2 c_0 \delta'_{0} + \sum_{m = 1}^\infty \frac{c_m}{\sqrt{\lambda_m}} \delta_{S_{\lambda_m}^1}\]
    satisfies
    \[\hat{\nu} \, = \, 2 \pi i \left(2(-1)^{k - 3} \, \hat{c}_0  \ts \alpha_k \ts \delta^{(k - 2)}_0 + \sum_{m = 1}^\infty \hat{c}_m \sum_{j = 0}^{\frac{k - 1}{2}} (-1)^{j - 1}  \beta_{j, k} \ts \lambda_m^{\frac{j - k + 1}{2}} \delta^{(j - 1)}_{S_{\hat{\lambda}_m}^1}\right).\]
    In particular, for $k = 3$ we obtain
    \[\hat{\nu} \, = \, i\left(-2\ts \hat{c}_0 \ts\delta'_0 + \sum_{m = 1}^\infty \frac{\hat{c}_m}{\sqrt{\lambda_m}} \, \delta_{S_{\hat{\lambda}_m}^1}\right).\]
    Thus, if $\hat{\mu} = i^{\varepsilon} \mu$, we get $\hat{\nu} = i^{\varepsilon + 1} \nu$.
\end{corollary}

\section{Recovering Poisson type summation formulas}\label{sec:poissontype}

\subsection{Modular forms for the Fricke group}

Let $N\geq 1$ and write
\[
\G_0(N)\, =\, \left\{\begin{pmatrix} a&b\\c&d\end{pmatrix} \in\mathrm{SL}_2(\ZZ)\, \mid \, c\equiv 0\pmod{N} \right\}.
\]
We also consider the Fricke group $\G_0^+(N)$, which is generated by $\G_0(N)$ and the Fricke involution
\[
    W_N \, =\, \begin{pmatrix}
        0& -1/\sqrt{N} \\\sqrt{N} &0
    \end{pmatrix}.
\]
For $M=\left(\begin{smallmatrix} a&b\\c&d\end{smallmatrix}\right)\in\G_0(4)$, we define the \emph{theta multiplier system} by 
\[
 v_\vartheta(M)\, =\, \left(\frac{c}{d}\right) \varepsilon_d^{-1},
 \]
where $\left(\frac{c}{d}\right)$ denotes the Kronecker symbol, and 
\[
 \varepsilon_d\, =\, \begin{cases}
     1,&\text{if }d\equiv 1 \pmod{4},\\
     i, & \text{if }d \equiv 3 \pmod{4}.
 \end{cases}
\]
It can be extended to another multiplier system of weight $1/2$ for $\G_0^+(4)$ that satisfies $\sqrt{-i \tau} = \sqrt{\tau} \ts v_\vartheta(W_4)$, which we will also denote by $v_\vartheta$.
\begin{example}
    The classical Poisson summation formula implies 
     that the Jacobi theta function $\vartheta(\tau) = \sum_{m \in \IZ} e(m^2 \tau)$
satisfies
\[ \vartheta(M \tau) \, = \, (c \tau + d)^{\frac{1}{2}} \ts v_{\vartheta}(M) \vartheta(\tau)\]
for all $M = \left(\begin{smallmatrix} a&b\\c&d\end{smallmatrix}\right)\in\G_0^+(4)$.
More generally, if $v_{\vartheta}^k$ is a multiplier system of weight $k / 2$ for $\Gamma_0^+(4)$, then
\[ \vartheta(\tau)^k \, = \, \sum_{m = 0}^\infty r_k(m) \ts e(m \tau) \] 
satisfies
\[ \vartheta^k(M \tau) \, =\,  (c \tau + d)^{\frac{k}{2}} \ts v_{\vartheta}^k(M) \ts \vartheta^k(\tau)\]
for all $M = \left(\begin{smallmatrix} a&b\\c&d\end{smallmatrix}\right)\in\G_0^+(4)$, where
    \begin{equation}\label{eq:representationnumbers}
        r_k(m) = \lvert \{ (x_1, \ldots, x_k) \in \IZ^k \ \vert\ x_1^2 + \ldots + x_k^2 = m \} \rvert.
    \end{equation}
\end{example}

The theta multiplier system is a special instance of more general multiplier systems for a subgroup $\Gamma \subseteq \SL_2(\IR)$ of weight $k$. These are maps $v : \Gamma \to \IC$ with $\lvert v \rvert = 1$ satisfying the cocycle relation
\[v(M_1 M_2) \, = \, \frac{j(M_1 M_2, \tau)^k}{j(M_1, M_2 \tau)^k \, j(M_2, \tau)^k} \ts v(M_1) v(M_2),\]
where $j\bigl(\left(\begin{smallmatrix}a & b \\ c & d\end{smallmatrix}\right), \tau\bigr) = (c \tau + d)$ (see Section 6 of \cite{iwaniec} for details). 

A second multiplier system of weight $1/2$ for $\SL_2(\IZ)$ is given by the \emph{eta multiplier system}, which is defined by
\[
    v_\eta(M)\, =\, {\begin{cases}e^{\frac {bi\pi }{12}}, &\mbox{if } c=0,\,d=1,\\[2pt]
    e^{i\pi \left({\frac {a+d}{12c}}-s(d,c)-{\frac {1}{4}}\right)}, & \mbox{if }c>0,\end{cases}}
\]
for $M=\left(\begin{smallmatrix}
    a&b\\c&d
\end{smallmatrix} \right)\in \Gamma_0^+(1) = \SL_2(\ZZ)$, where $s(d,c)$ is the Dedekind sum
\[
s(d,c)\, =\, \sum _{n=1}^{c-1}{\frac {n}{c}}\left({\frac {dn}{c}}-\left\lfloor {\frac {dn}{c}}\right\rfloor -{\frac {1}{2}}\right).
\]

\begin{definition}\label{def:modformthetamult}
    Let $k \in \frac{1}{2} \IZ$ and $v : \Gamma_0^+(N) \to \IC$ be a multiplier system of weight $k$ for $\Gamma_0^+(N)$. A holomorphic function $f : \IH \to \IC$ is called \emph{modular form of weight $k$ for $\Gamma_0^+(N)$ and multiplier system $v$} if the following hold:
    \begin{enumerate}
 \item For all $M=\left(\begin{smallmatrix} a&b\\c&d\end{smallmatrix}\right)\in\G_0^+(N)$, one has
 \[
  f\left(M \tau \right) \, =\, (c\tau+d)^{k} v(M) f(\tau).
 \]
 \item $f$ is holomorphic at the cusps of $\G_0^+(N)$.
\end{enumerate}
\end{definition}

A modular form $f$ of weight $k$ for $\Gamma_0^+(N)$ and multiplier system $v$ has a Fourier expansion of the form
\[f(\tau) \, =\,  \sum_{m \in \IZ}^\infty c_m \ts e\bigl((m + m_0) \tau\bigr),\]
where $v(T) = e(-m_0), \ 0 \leq m_0 < 1$.
This can be seen by applying the transformation $T = \left(\begin{smallmatrix} 1&1\\0&1\end{smallmatrix}\right)$. The second condition of Definition \ref{def:modformthetamult} ensures that $c_m = 0$ for $m < 0$. Moreover, the Fourier coefficients $c_m$ have polynomial growth.

Note that the complex vector space of modular forms of fixed weight, fixed group, and given multiplier system is finite-dimensional. 
For instance, we have the following dimension formula for modular forms of half-integral weight $k$ for $\Gamma_0(4)$ and the theta multiplier (compare \cite[Proposition 15.1.2]{cohenstroemberg})
\begin{align}\label{eq:dimmodforms}
    \mathrm{dim}(M_k(\Gamma_0(4))&=\begin{cases}
        0, &\text{when }k<0,\\
        1+\lfloor\frac{k}{2}\rfloor, & \text{when } k\geq 0.
    \end{cases}
\end{align}
For dimension formulas of the Fricke group, see \cite{YangZhanFrickeDimensions} and references therein.

\begin{remark}
We note that dimension formulas for the spaces of modular forms as above can be used to determine the dimension of certain subspaces of the space of $k$-spherical measures. For example, the space of modular forms of weight $1/2$ for $\widetilde\G_0^+(4)$ and the theta multiplier is one-dimensional by \eqref{eq:dimmodforms}. This implies that the only crystalline measure supported on the set $\{0,\pm\sqrt{m}\,|\, m\in \NN\}$ satisfying $\mu=\hat\mu$ is the standard Dirac comb $\delta_{\ZZ}$. 
\end{remark}

\subsection{Poisson type summation formulas involving the coefficients of modular forms}
In this section, we review several results in the literature and explain how these follow from our results.

\subsubsection{Results of Cohn--Gon\c{c}alves and Gon\c{c}alves}
\begin{proposition}\label{lem:ModularFormsAndRadialEigenmeasures}
    Let $f : \IH \to \IC, \, f(\tau) = \sum_{m = 0}^\infty c_m \ts e(m \tau / b)$, be a modular form of weight $\frac{k}{2}$ for the Fricke group $\Gamma_0^+(N)$ with multiplier system $v$. Then the measure
    \[\mu \, = \, \sum_{m = 0}^\infty \frac{c_m}{\vol\bigl(S_{\sfrac{2 m}{(b\sqrt{N}})}^{^{\ts k}}\bigr)} \ts \delta_{S_{\sfrac{2 m}{(b\sqrt{N}})}^{^{\ts k}}}\]
    satisfies $\hat{\mu} = \frac{v(W_N)}{v_{\vartheta}^k(W_4)} \mu$. In particular, for all $\varphi \in \calS(\IR^k)^{O(k)}$, we have the summation formula
    \[\sum_{m = 0}^\infty c_m \ts \varphi\!\left(\!\sqrt{\frac{2 m }{ b \sqrt{N}}}\ts \right) \, = \,  \frac{v(W_N)}{v^k_{\vartheta}(W_4)} \sum_{m = 0}^\infty c_m \ts \hat{\varphi}\left(\!\sqrt{\frac{2 m }{ b \sqrt{N}}}\right).\]
\end{proposition}

\begin{proof}
    Consider $\tilde{f}(\tau) = f\bigl(\tau / \sqrt{N} \ts \bigr)$. Then
    \[\tilde{f}(\tau) \, = \, \sum_{m = 0}^\infty c_m \ts e\!\left(\frac{m \tau}{ b\sqrt{N}}\right)\]
    and $\tilde{f}(-1 / \tau) = \tfrac{v(W_N)}{v^k_{\vartheta}(W_4)} \sqrt{-i \tau}^{\, k} \tilde{f}(\tau)$. Applying Theorem \ref{lem:IsoBetweenMeasuresAndFunctions} to $\mu_{\tilde{f}}$ yields $\hat{\mu} = \frac{v(W_N)}{v_\vartheta^k(W_4)}\mu$. For the summation formula use that, for $\varphi \in \calS(\IR^k)^{O(k)}$, we have $\widehat{\hat{\varphi}\ts}(x) = \varphi(-x) = \varphi(x)$, which implies
    \[ \langle \mu, \varphi \rangle \, = \,  \langle \mu, \hat{\hat{\varphi\,}} \ts \rangle \, = \, \langle \hat{\mu}, \hat{\varphi} \ts\rangle. \qedhere\]
\end{proof}

\begin{remark}
 Applying Proposition \ref{lem:ModularFormsAndRadialEigenmeasures} to the Jacobi theta function $\vartheta$ yields the usual Poisson summation formula. More generally, we obtain the summation formulas
\[\sum_{m = 0}^\infty r_k(m) \ts \varphi(m) \, = \, \sum_{m = 0}^\infty r_k(m) \ts \hat{\varphi}(m)\]
for all $\varphi \in \calS(\IR^k)^{O(k)}$.
\end{remark}
\begin{remark}
We apply Proposition \ref{lem:ModularFormsAndRadialEigenmeasures} to the weight $6$ Eisenstein series $E_6$ given by
    \[
       E_6(\tau) \, = \, 1-504\sum_{m=1}^\infty \sigma_5(m)\ts e(m\tau),
    \]
    where $\sigma_k(m)=\sum_{d\mid m} d^{k}$ is the divisor sum. This recovers Lemma 2.1 of \cite{cohngonzalves} as $E_6(\tau)$ is a modular form of weight $6$ for $\SL_2(\ZZ)$ with a trivial multiplier system.
 
\end{remark}
\begin{remark}
    The construction of crystalline measures of \cite[Lemma 19]{goncalves} is the special case $k = 1$ of Theorem \ref{lem:IsoBetweenMeasuresAndFunctions}. In particular, the crystalline measures constructed in \cite{goncalves} using eta products of weight $1/2$ are recovered by applying Proposition \ref{lem:ModularFormsAndRadialEigenmeasures} to the corresponding eta products in \cite[Appendix]{goncalves}. We note that Proposition \ref{lem:ModularFormsAndRadialEigenmeasures} works for arbitrary half-integral weights (see \cite[Section 2.2]{koehler} for an infinite family of eta products satisfying the above assumptions).
\end{remark}
\begin{remark}
 Let $\chi : \IZ / N \IZ \to \IC$ be a Dirichlet character.  For $\gamma = \left(\begin{smallmatrix}
        a & b \\ c & d
    \end{smallmatrix}\right) \in \Gamma_0(N)$, set $\chi(\gamma) = \chi(d)$. This defines a multiplier system of $\Gamma_0(N)$ for every integral weight. We call $\chi$ quadratic if $\chi^2$ is the trivial character. Since $\Gamma_0(N)$ has index $2$ in $\Gamma_0^+(N)$, we can extend $\chi$ to a character of $\Gamma_0^+(N)$ in exactly two ways in this case. If $\chi(-\id) = 1$, we have the possible choices $\chi_{\pm}(W_N) = \pm1$; if $\chi(-\id) = -1$, we can choose between $\chi_{\pm}(W_N) = \pm i$. Both define multiplier systems of $\Gamma_0^+(N)$ for every integral weight.

    Let now $N$ be an odd prime and $\chi : \IZ / N \IZ \to \IC$ be the Dirichlet character given by the Legendre symbol $\chi(d) = \left(\frac{d}{N}\right)$. Let $k$ be even such that $k / 2 \geq 3$ and assume $\chi(-1) = (-1)^{\frac{k}{2}}$. Let
    \[A_{\frac{k}{2}}(N) \, = \, (-1)^{\left \lfloor \frac{k}{4} \right \rfloor} \frac{N^{\frac{k - 1}{2}}(k / 2 - 1)!}{(2 \pi)^{\frac{k}{2}}} L(\chi, \, k / 2), \quad \text{where }\ L(\chi, \, k / 2) = \sum_{m = 1}^\infty \chi(m)\, m^{-\frac{k}{2}}.\]
    By \cite[Proposition 1.10]{koehler}, the function
    \[E^{}_{\frac{k}{2},\, P, \, \pm i}\ts(\tau) \, = \, 1 + \frac{1}{A_{\frac{k}{2}}(N)} \sum_{m = 1}^\infty \bigg(\sum_{d \mid m} \bigl(\chi(d) \pm N^{\frac{k/2 - 1}{2}} \chi(n / d)\bigr)\, d^{\, k/2 - 1}\bigg)\ts e(m \tau)\]
    is a modular form of weight $k / 2$ for $\Gamma_0^+(N)$, and multiplier system given by $\chi$ on $\Gamma_0(N)$ and extended by $\chi(W_N) = \mp (-i)^{k/2} (-1)^{\left\lfloor \frac{k}{4}\right\rfloor}$ to all of $\Gamma_0^+(N)$. Hence we obtain Fourier eigenmeasures with eigenvalue $\frac{\chi(W_N)}{v_\vartheta^{k}(W_4)}$.
\end{remark}

\subsubsection{Results of Lev--Reti and Meyer}

\begin{proposition}\label{lem:ModularFormsAnd1DimEvenMeasures}
Let $k \geq 1$ be odd and $f : \IH \to \IC, f(\tau) = \sum_{m = 0}^\infty c_m e\bigl(\tfrac{m \tau}{b}\bigr)$ be a modular form of weight $\frac{k}{2}$ for the Fricke group $\Gamma_0^+(N)$ with multiplier system $v$. Then, for even $\varphi \in \calS(\IR)$, we have the summation formula
\begin{align*}
    &c_0 \, \varphi(0) + \sum_{m = 1}^\infty c_m \, \varphi\left(\!\sqrt{\frac{2 m }{ b \sqrt{N}}}\ts\right) \\
    &= \frac{v(W_N)}{v^k_{\vartheta}(W_4)} \left(\alpha_k \ts c_0 \ts  \hat{\varphi}^{(k - 1)}(0) + \sum_{m = 1}^\infty c_m \sum_{j = 0}^{\frac{k - 1}{2}} \beta_{j, k} \left(\frac{2 m}{b \sqrt{N}}\right)^{\frac{j - k + 1}{2}} \hat{\varphi}^{(j)}\left(\!\sqrt{\frac{2 m }{ b \sqrt{N}}}\ts\right)\right).
\end{align*}
\end{proposition}

\begin{proof}
    We consider $\tilde{f}(\tau) = f\bigl(\tau / \sqrt{N}\,\bigr)$, and apply Lemma \ref{lem:EvenMeasureFromSphericalMeasure} to $\mu_{\tilde{f}}$. The summation formula is obtained using $\hat{\hat{\varphi}\,} = \varphi$ for even $\varphi \in \calS(\IR)$.
\end{proof}

\begin{proposition}\label{lem:ModularFormsAnd1DimOddMeasures}
Let $k \geq 3$ be odd and $f : \IH \to \IC$, $f(\tau) = \sum_{m = 0}^\infty c_m e\bigl(\tfrac{m\tau}{b}\bigr)$ be a modular form of weight $\frac{k}{2}$ for the Fricke group $\Gamma_0^+(N)$ with multiplier system $v$. Then, for odd $\varphi \in \calS(\IR)$, we have the summation formula
\begin{align*}
    &c_0\ts \varphi'(0) + \sum_{m = 1}^\infty \frac{c_m}{\sqrt{2 m / (b \sqrt{N})}}\, \varphi\!\left(\!\sqrt{\frac{2 m }{ b \sqrt{N}}}\ts\right) \\
    &= 2 \pi i\,\frac{v(W_N)}{v^k_{\vartheta}(W_4)} \!\left(c_0 \ts \alpha_k \ts\hat{\varphi}^{(k - 2)}(0) + \sum_{m = 1}^\infty c_m \sum_{j = 0}^{\frac{k - 1}{2}} \beta_{j, k} \left(\frac{2 m}{b \sqrt{N}}\right)^{\frac{j - k + 1}{2}} \hat{\varphi}^{(j - 1)}\left(\!\sqrt{\frac{2 m }{ b \sqrt{N}}}\ts\right)\right).
\end{align*}
\end{proposition}

\begin{proof}
Again, consider $\tilde{f}(\tau) = f\bigl(\tau / \sqrt{N} \,\bigr)$ and apply Corollary \ref{cor:OddMeasureFromSphericalMeasure} to $\mu_{\tilde{f}}$. The summation formula is then obtained by using $\hat{\hat{\varphi}\,} = - \varphi$ for odd $\varphi \in \calS(\IR)$.
\end{proof}

\begin{remark}\label{rem:thetak}
    We apply Proposition \ref{lem:ModularFormsAnd1DimOddMeasures} to the functions $\vartheta^k(\tau)$. This gives the following summation formula, which is the main result of \cite{levreti}, and generalizes a theorem by Meyer \cite{meyerpnas},
    \begin{align*}
    &\varphi'(0) + \sum_{m = 1}^\infty \frac{r_k(m)}{\sqrt{m}} \,\varphi(\sqrt{m}\ts) \\
    &\quad\quad= 2 \pi i \left(\alpha_k \ts \hat{\varphi}^{(k - 2)}(0) + \sum_{m = 1}^\infty r_k(m) \sum_{j = 0}^{\frac{k - 1}{2}} \beta_{j, k} \, m^{\frac{j - k + 1}{2}} \hat{\varphi}^{(j - 1)}(\sqrt{m}\ts)\right).
\end{align*}
    In particular, for $k = 3$, we recover Guinand's example as in Section \ref{subsec:guinand}.
\end{remark}

\section{Measures and Fourier series in higher dimensions}\label{sec:measuresfourierhigher}
Now, we generalize the results of Section \ref{sec:MeasuresFromFourierSeries} to higher dimensions. To do so, let $\underline{k} = (k_1, \ldots, k_n) \in \IN^n, n \in \IN$ and let $\IR^{\underline{k}} = \prod_{i = 1}^n \IR^{k_i}$. We obtain a representation of $\Mp_{2}(\IR)^n$ on $\calS(\IR^{\underline{k}})$. For $\underline{\tau} = (\tau_1, \ldots, \tau_n) \in \IH^n$, consider the Gaussian
\[g(\underline{x}, \ts \underline{\tau}) \, = \, \prod_{i = 1}^n g(x_i, \tau_i) \in \calS(\IR^{\underline{k}})^{O(\underline{k})},\]
where $O(\underline{k}) = O(k_1) \times \ldots \times O(k_n)$.

We have the following generalization of the case $n = 1$. Its proof is analogous to the one of Lemma \ref{lem:Mp2Equivariance}, again, using the denseness of the Gaussians in the space of Schwartz functions that are radial in every variable \cite[Proposition 4.1]{RadchenkoStoller}.

\begin{lemma}
For $\mu \in \calS'(\IR^{\underline{k}})$ and $M \in \Mp_2(\IR)^n$, we have $\theta_{\mu \vert M} = \theta_{\mu}\vert_{\underline{k} / 2} M$, i.e.\@, the map $\mu \mapsto \theta_\mu = \langle \mu, g(\cdot, \cdot)\rangle$ is $\Mp_2(\IR)^n$-equivariant. Moreover, if $\mu \in \calS'(\IR^{\underline{k}})^{O(\underline{k})}$, then $\hat{\mu} = \mu$ if and only if
\[\theta_\mu(-1/ \underline{\tau}) \, = \, \sqrt{-i \underline{\tau}}^{\,\underline{k}} \,\theta_\mu(\underline{\tau}).\]
\end{lemma}

For $\underline{\lambda} = (\lambda_1, \ldots, \lambda_n) \in \IR_{\geq 0}^n$, define $S_{\underline{\lambda}}^{\underline{k}} = \prod_{i = 1}^n S_{\lambda_i}^{k_i}$. Let $\Meas_{\underline{k}}$ be the space of measures $\mu$ on $\IR^{\underline{k}}$ such that there exist locally finite sequences $\lambda_{\underline{m}}, \ \hat{\lambda}_{\underline{m}} \in \IR_{\geq 0}^n$ and at most polynomially growing sequences $c_{\underline{m}}, \ \hat{c}_{\underline{m}} \in \IC$ with
\[\mu = \sum_{\underline{m} = 0}^\infty c_{\underline{m}} \ts \delta_{S_{\lambda_{\underline{m}}}^{^{\ts \underline{k}}}} \qquad \text{and} \qquad \hat{\mu} = \sum_{\underline{m} = 0}^\infty \hat{c}_{\underline{m}} \ts \delta_{S_{\hat{\lambda}_{\underline{m}}}^{^{\ts \underline{k}}}}.\]
We call these measures \emph{$\underline{k}$-spherical measures}. Similarly, we define the space $\FS_{\underline{k}}$. Then, we obtain the following analogue of the $n = 1$ case.

\begin{theorem}\label{thm:mainhigherdim}
The map
\[\theta : \, \Meas_{\underline{k}} \, \to \,  \FS_{\underline{k}}\ts, \quad \mu \mapsto \theta_\mu(\underline{\tau}) = \langle \mu, g(\cdot, \underline{\tau}) \rangle\]
is an $\Mp_2(\IR)^n$-equivariant isomorphism. The inverse is given by
\[f(\underline{\tau}) \, = \,  \sum_{\underline{m} = 0}^\infty c_{\underline{m}} \ts e(\lambda_{\underline{m}} \ts \underline{\tau} / 2) \ \longmapsto \ \mu_f = \sum_{\underline{m} = 0}^\infty \, \frac{c_{\underline{m}}}{\vol\bigl(S_{\lambda_{\underline{m}}}^{^{\ts\underline{k}}}\bigr)} \delta_{S_{\lambda_{\underline{m}}}^{^{\ts \underline{k}}}}.\]
Under this isomorphism, measures $\mu \in \Meas_{\underline{k}}$ with $\hat{\mu} = \mu$ correspond to functions $f \in \FS_{\underline{k}}$ with $f(-1/\underline{\tau}) = \sqrt{-i \underline{\tau}}^{\, \underline{k}} f(\underline{\tau})$.
\end{theorem}

Let $i_* : \calS(\IR^n) \to \calS(\IR^{\underline{k}})^{O(\underline{k})}$ be the map given by $(i_* \varphi)(x_1, \ldots, x_n) = \varphi(\lvert x_1 \rvert, \ldots, \lvert x_n \rvert)$. We have

\begin{lemma}[{\cite[Theorem 1.1]{GrafakosFourierMultiradial}, \cite[Remark 4.3]{levreti}}]\label{lem:fouriertrafohigherdim}
Let $\underline{k}$ be odd, i.e.\@, $k_i$ is odd for all $1 \leq i \leq n$, and let $\varphi \in \calS(\IR^n)^{O(\underline{1})}$. Then the Fourier transform of $i_* \varphi$ is given by
\[\widehat{i_* \varphi \,}(\underline{x}) \, = \, \sum_{\underline{j} = 0}^{\frac{\underline{k} - \underline{1}}{2}} \gamma_{\underline{j}, \, \underline{k}}(x) \ts i_* \bigl(D_{\underline{j}, \, \underline{k}}(\underline{x}) \hat{\varphi}\ts \bigr)(\underline{x}),\]
where
\[\gamma_{\underline{j}, \underline{k}}(\underline{x}) = \prod_{\substack{i = 1}}^n \gamma_{j_i, \, k_i}(x_i), \qquad D_{\underline{j}, \underline{k}}(\underline{x}) = \prod_{\substack{i = 1}}^n D_{j_i, \, k_i}(x_i).\]
\end{lemma}

\begin{proof}
This is \cite[Theorem 1.1]{GrafakosFourierMultiradial} together with the observation of \cite[Remark 4.3]{levreti}, or an application of Lemma \ref{lem:FourierTransformRadial}.
\end{proof}

If $\mu$ is an $O(\underline{k})$-invariant distribution on $\calS(\IR^{\underline{k}})$, we define a tempered distribution $i^* \mu$ on $\calS(\IR^n)$ via $\langle i^* \mu, \varphi \rangle = \langle \mu, i_* \varphi \rangle$. We obtain the following analogues of Lemma \ref{lem:EvenMeasureFromSphericalMeasure} and Lemma \ref{lem:OddMeasureFromSphericalMeasure}.

\begin{corollary}\label{cor:higher1}
Let $\underline{k}$ be odd and $\mu \in \calS'(\IR^{\underline{k}})^{O(\underline{k})}$. Then
\[\widehat{i^* \mu \,} \, = \, \sum_{\underline{j} = 0}^{\frac{\underline{k} - 1}{2}} (-1)^{\underline{j}} \ts D_{\underline{j}, \,\underline{k}}(\underline{x}) \ts \gamma_{\underline{j}, \, \underline{k}}(\underline{x}) \ts i^* \hat{\mu}.\]
\end{corollary}

\begin{corollary}\label{cor:higher2}
Let $\underline{k}$ be odd and $\mu \in \calS'(\IR^{\underline{k}})^{O(\underline{k})}$. Let $\Delta \subseteq \{1, \ldots, n\}$ and further define $x_\Delta = \prod_{i \in \Delta} x_i$. Then
\[\widehat{\frac{1}{x_\Delta} i^* \mu} \, =\,  i \sum_{\underline{j} = 0}^{\frac{\underline{k} - 1}{2}} (-1)^{j+\lvert \Delta \rvert} \tilde{D}_{\underline{j},\, \underline{k}, \,\Delta}(\underline{x}) \gamma_{\underline{j},\, \underline{k}}(\underline{x}) i^* \hat{\mu},\]
where
\[\tilde{D}_{\underline{j}, \,\underline{k},\, \Delta}(\underline{x}) = \prod_{\substack{i = 1 \\ i \in \Delta}}^n \tilde{D}_{j_i, \,k_i}(x_i) \prod_{\substack{i = 1 \\ i \notin \Delta}}^n D_{j_i, \, k_i}(x_i).\]
\end{corollary}

\subsection{An example using Hilbert modular forms}
\subsubsection{Hilbert modular forms}
Denote by $F$ a totally real number field of degree $n$ over~$\mathbb{Q}$. Let $\mathcal{O}_F$ be its ring of integers and $\{\sigma_1,\ldots,\sigma_n\}$ its embeddings. The group $\Gamma_F=\mathrm{SL}_2(\mathcal{O}_F)$ acts on 
\[
\mathbb{H}^n \, = \, \{\underline{\tau}=(\tau_1,\ldots,\tau_n)\, \mid \,\tau_j=u_j+i\ts v_j,\,v_j>0 \text{ for all } 1\leq j\leq n\}
\]
properly discontinuously via 
\[
\begin{pmatrix}
a&b\\c&d
\end{pmatrix} (\tau_1,\ldots,\tau_n) \, =\, \left(\frac{\sigma_1(a) \tau_1+\sigma_1(b)}{\sigma_1(c)\tau_1+\sigma_1(d)},\, \ldots, \, \frac{\sigma_m(a) \tau_1+\sigma_m(b)}{\sigma_n(c)\tau_1+\sigma_n(d)}\right).
\]
The quotient $Y_F=\Gamma_F\setminus\mathbb{H}^n$ is a complex variety that can be compactified (Bailey--Borel compactification). 

A holomorphic function $f: \mathbb{H}^n\to \CC$ is called a \textit{holomorphic Hilbert modular form of weight $\underline{k}=(k_1,\ldots,k_n)\in\ZZ^n$ for $\Gamma_F$} if 
\[
 f(\underline{\tau}) \,=\,  j_{\underline{k}}^{-1} (\gamma,\underline{\tau})f(\gamma \underline{\tau}),
\]
where $ j_{\underline{k}}(\gamma,\underline{\tau})= \prod_{j=1}^n \bigl(\sigma_j(c)\tau_j+\sigma_j(d)\bigr)^{k_j}$. Hilbert modular forms also possess Fourier expansions in the cusps. Similarly, one can define Hilbert modular forms with multiplier systems.

\begin{example}
    Let $F$ be a real quadratic field. Let $\frakd_F$ be the different ideal of $F$ and $D_F$ be the discriminant of $F$. Define the Dedekind zeta function by
    \[\zeta_F(s) \, = \, \sum_{\fraka \subseteq \calO_F} N(\fraka)^{-1},\]
    where $N(\fraka) = \lvert \calO_F / \fraka \rvert$. It has a meromorphic continuation to all $s \in \IC$ with a simple pole at $s = 1$. Then, for $k \geq 2$ even, the Eisenstein series
    \[E_k(\underline{\tau}) \, = \, \zeta(1 - k) + \frac{(2 \pi i)^{2k}}{(k - 1)!^2}\ts  D_F^{^{\frac{1}{2} - k}} \sum_{\substack{m \in \frakd_F^{-1} \\ m \gg 0}} \sigma_{k - 1}(\frakd_F m) \ts e(\tr(\underline{\tau} m^2))\]
    is a Hilbert modular form of weight $\underline{k} = (k, k)$, where
    \[\sigma_{s}(\frakm) \, = \, \sum_{\frakn \mid \frakm} N(\frakn)^s,\]
    see for example \cite[Section 1.5]{Bruinier123}. For $m \in \calO_F$, write $\underline{m} = (\sigma_1(m), \sigma_2(m))$. Then the measure
    \[\mu \, = \, \sum_{\substack{m\in \frakd_F^{-1} \\ m \gg 0}} \frac{\sigma_{k - 1}(\frakd_F \ts m)}{\vol\bigl(S_{2\underline{m}}^{^{\,\underline{k}}}\bigr)} \delta_{S_{2\underline{m}}^{^{\,\underline{k}}}}\]
    is a Fourier eigenmeasure.
\end{example}

\section{Further constructions}\label{sec:further}

The usual operations on the function spaces $\FS_{\underline{k}}$ induce operations on $\Meas_{\underline{k}}$.

\begin{lemma}
    Let $\underline{k} \in \IN^n, \ \underline{k}' \in \IN^{n'}$. The map
    \[\FS_{\underline{k}} \times \FS_{\underline{k}'}\, \to \, \FS_{(\underline{k},\, \underline{k}')}, \qquad (f, f'\ts ) \,\mapsto \, \bigl(f \otimes f' : (\tau, \tau') \mapsto f(\tau) f(\tau')\bigr)\]
    induces a map
    \[\Meas_{\underline{k}} \times \Meas_{\underline{k}'} \, \to \, \Meas_{(\underline{k}, \,\underline{k}')}, \quad (\mu, \mu) \, \mapsto \, \mu \otimes \mu'\]
    such that $\theta_{\mu \otimes \mu'} = \theta_\mu \otimes \theta_{\mu'}$. Explicitly, the map is given by
    \[\left(\sum_{\underline{m} = 0}^\infty c_{\underline{m}}\, \delta_{S_{\lambda_{\underline{m}}}^{^{\ts \underline{k}}}}, \, \sum_{\underline{m}' = 0}^\infty c_{\underline{m}'}'\, \delta_{S_{\lambda_{\underline{m}'}'}^{^{\ts \underline{k}'}}}\right) \ \longmapsto \ \sum_{\underline{m} = 0}^\infty \sum_{\underline{m}' = 0}^\infty c_{\underline{m}} \ts c_{\underline{m}'}'\ts  \delta_{S_{\lambda_{\underline{m}}}^{^{\ts \underline{k}}}  \times \, S_{\lambda_{\underline{m}'}'}^{^{\ts \underline{k}'}}}.\]
    Both maps respect the action of $\Mp_2(\IR)^n \times \Mp_2(\IR)^{n'} \simeq \Mp_2(\IR)^{n + n'}$.
\end{lemma}

\begin{lemma}
    Given a partition $\Delta = \{\Delta_1, \, \ldots, \, \Delta_{n'}\}$ of $\{1, \ldots, n\}$ with $n' = \lvert \Delta \rvert$, there is a~map $\diag_\Delta : \IH^{n'} \to \IH^n$ which induces a map
    \[\diag_\Delta^* : \, \FS_{\underline{k}} \, \to \, \FS_{\underline{k}'},\]
    where $\underline{k}' = \bigl(\sum_{j \in \Delta_i} k_j \bigr)_{1 \leq i \leq n'} \in \IN^{n'}$. There also exists a map
    \[\diag_\Delta^* :\, \Meas_{\underline{k}} \, \to \, \Meas_{\underline{k}'},\]
    such that
    \[\diag_\Delta^*(\theta_{\mu}) \, = \, \theta_{\diag_\Delta^*(\mu)}.\]
    Explicitly, the map is given by
    \[\sum_{\underline{m} = 0}^\infty c_{\underline{m}} \ts \delta_{S_{\lambda_{\underline{m}}}^{^{\ts \underline{k}}}} \ \longmapsto \ \sum_{\underline{m} = 0}^\infty c_{\underline{m}} \,\frac{\vol\bigl({S_{\lambda_{\underline{m}}}^{^{\ts \underline{k}}}}\bigr)}{\vol\bigl(S_{\lambda_{\underline{m}}}^{^{\ts \underline{k}'}} \bigr)}  \,\delta_{S_{\lambda_{\underline{m}}}^{^{\ts \underline{k}'}}}.\]
    If we embed $\Mp_2(\IR)^{n'}$ into $\Mp_2(\IR)^n$ appropriately, then all these maps respect the action of $\Mp_2(\IR)^{n'}$.
\end{lemma}

\noindent Combining these two maps, we have

\begin{lemma}
    For $\underline{k}, \, \underline{k}' \in \IN^n$, we have a map
    \[\Meas_{\underline{k}} \times \Meas_{\underline{k}'}\, \to \, \Meas_{\underline{k} + \underline{k}'}\]
    which is induced by the map
    \[\FS_{\underline{k}} \times \FS_{\underline{k}'} \, \to \, \FS_{\underline{k} + \underline{k}}\, , \qquad (f, f'\ts ) \mapsto \bigl(\tau \mapsto f(\tau) f'(\tau)\bigr),\]
    and respects the action of $\Mp_2(\IR)^n$ if we embed it diagonally into $\Mp_2(\IR)^n \times \Mp_2(\IR)^n$.
\end{lemma}

\begin{example}
    Consider the function $f(\tau, \tau'\ts ) = \vartheta(\tau / 2) \ts \vartheta(\tau' / 2) \in \FS_{(1/2, \ts 1/2)}$. The corresponding measure is the usual Dirac comb $\delta_{\IZ^2} = \delta_{\IZ} \otimes \delta_{\IZ}$, and we have $\theta_{\delta_{\IZ^2}}(\tau, \tau'\ts ) = f(\tau, \tau'\ts )$. Using the diagonal embedding $\diag : \IH \to \IH^2, \, \tau \mapsto (\tau, \tau)$, we have
    \[(\diag^* f)(\tau) \, = \, \vartheta(\tau / 2)^2 \, = \, \theta_{\delta_{\IZ}}(\tau)^2 = \theta_{\diag^*(\delta_{\IZ^2})}.\]
    
\end{example}

\begin{lemma}
    Let $\underline{k} \in \IN^n$ and $k = \sum_{j = 1}^n k_j$. Then there is a diagonal embedding of $\Mp_2(\IR)$ into $\Mp_2(\IR)^n$ and the actions of the orthogonal group $O(k)$ on $\calS(\IR^k), \calS'(\IR^k)$ commute with the actions of the diagonally embedded $\Mp_2(\IR)$. In particular, if $\mu \in \calS'(\IR^k)$ is a Fourier eigendistribution, then for every $\gamma \in O(k)$, the tempered distribution $\gamma \mu$ is a Fourier eigendistribution.
\end{lemma}

\begin{example}
    Let $\gamma \in O(2)$ be a rotation through $-\alpha$. Then $\gamma \delta_{\IZ^2} = \delta_{\gamma^{-1} \IZ^2}$. Assume that $\tan(\alpha) = \lambda' = - 1 / \lambda$ is a root of a polynomial $X^2 + A X - 1$ and let $\lambda$ be its second root. Let $K = 1 / \sqrt{1 + \lambda^2},\, K' = 1 / \sqrt{1 + {\lambda'}^2}$. Then, the corresponding theta function is given by
    \[\theta_{\gamma^{-1} \delta_{\IZ^2}}(\tau, \tau') \, = \, \sum_{x \in \IZ[\lambda]} e\bigl(K^2 x^2 \tau / 2 + {K'\ts}^2 {x'\ts}^2 \tau' / 2\bigr),\]
    where $x' = a + b \lambda'$ if $x = a + b \lambda$. By the previous lemma, we have $\theta_{\gamma^{-1} \delta_{\IZ^2}}(-1/\tau,\ts -1/\tau') = \sqrt{-i \tau \tau'}\,\theta_{\gamma^{-1} \delta_{\IZ^2}}(\tau, \tau')$. In fact, it is a Hilbert modular form of weight $(1/2,\, 1/2)$ for the number field $\IQ(\lambda)$ and equal to a Hilbert theta function for the lattice $\IZ[\lambda]$.
\end{example}

As the last construction, we mention that, if $\underline{k} \in \IN^n,\,  \underline{k}' \in \IN^{n'}$ and $f \in \FS_{(\underline{k}, \underline{k'})}$, then $f'(\ts\underline{\tau}'\ts) = f(i, \ldots, i, \,\underline{\tau}'\ts)$ satisfies $\sqrt{-i \underline{\tau}'}^{\, -\underline{k}'} f'(-1/\underline{\tau}'\ts ) = \tilde{f}'(\ts \underline{\tau}'\ts )$, where $\tilde{f}'(\ts \underline{\tau}) = f(i, \ldots, i, \, \underline{\tau}'\ts )$. Note that $f(i, \ldots, i, \,\underline{\tau}'\ts)$ is usually not in $\Meas_{\underline{k}'}$, since the support of the Fourier coefficients will not be locally finite. On the side of $\mu \in \Meas_{(k,\, \underline{k}')}$, this corresponds to considering the measure $\nu$ on $\IR^{\underline{k}'}$ given by
\[\langle \nu, \varphi \rangle \, = \, \langle \mu, \ts g\bigl(x_{\underline{k}}, \ts (i, \ldots, i)\bigr)\ts  \varphi(x_{\underline{k}'})\rangle\]
for $\varphi \in \calS(\IR^{\underline{k}'})$. This approach can be used to construct eigenmeasures, see \cite{baakeeigenmeasure}, or when studying one-dimensional Dirac combs with dense support, confer \cite{richard}.

\section*{Acknowledgements}
The work on this article started with a letter from Yves Meyer to Michael Baake in which Meyer observed the connection between certain crystalline measures and modular forms of weight $1/2$. We thank both of them for sharing their inspiration and insight. Moreover, we thank Michael Baake for many inspiring and helpful conversations.

\bibliographystyle{alpha}
\bibliography{bib.bib}
\end{document}